\documentclass[final,onefignum,onetabnum]{siamart171218}
\pdfoutput=1
\usepackage{amsmath}
\usepackage{amsfonts}
\usepackage{bm}
\usepackage{lipsum}
\usepackage{enumerate}
\usepackage{graphicx}
\usepackage{amsfonts}
\usepackage{mathrsfs}
\usepackage{subfigure}
\usepackage{epstopdf}
\usepackage{url}
\usepackage{verbatim}
\usepackage{amssymb}
\usepackage{hyperref}
\usepackage{algorithmic}
\usepackage{tikz}
\usetikzlibrary{shapes,calc,arrows,arrows.meta,decorations.markings}
\usepackage{color}
\usepackage{listings}
\usepackage[framemethod=tikz]{mdframed}

\newcommand{\tx}{\tilde{x}}

\title{Quadrature by Two Expansions: Evaluating Laplace Layer Potentials using Complex Polynomial and 
Plane Wave Expansions}
\author{
  Lingyun Ding\thanks{Department of Mathematics, the University of North Carolina at Chapel Hill (\email{dingly@live.unc.edu},
   \email{huang@email.unc.edu}, \email{marzuola@email.unc.edu}).}
 \and Jingfang Huang\footnotemark[1] \thanks{Corresponding author.}
 \and Jeremy L. Marzuola\footnotemark[1]
 }

\begin{document}
\maketitle
\tableofcontents

\begin{abstract}
The recently developed quadrature by expansion (QBX) technique \cite{klockner2013quadrature} 
accurately evaluates the layer potentials with singular, weakly or nearly singular, or even 
hyper singular kernels in the integral equation reformulations of partial differential equations.
The idea is to form a local complex polynomial or partial wave expansion centered at a 
point away from the boundary to avoid the singularity in the integrand, 
and then extrapolate the expansion at points near or even exactly on the boundary. 
In this paper, in addition to the local complex Taylor polynomial expansion, we derive new 
representations of the Laplace layer potentials using both the local complex polynomial and 
plane wave expansions. Unlike in the QBX, the local complex polynomial expansion in the 
new quadrature by two expansions (QB2X) method only collects the far-field contributions and 
its number of expansion terms can be analyzed using tools 
from the classical fast multipole method. The plane wave type expansion in the QB2X method
better captures the layer potential features near the boundary. It is derived 
by applying the Fourier extension technique to the density and boundary geometry functions 
and then analytically utilizing the Residue Theorem for complex contour integrals. The
internal connections of the layer potential with its density function and curvature on 
the boundary are explicitly revealed in the plane wave expansion and its error is bounded 
by the Fourier extension errors. We present preliminary numerical results to demonstrate 
the accuracy of the QB2X representations and to validate our analysis.
\end{abstract}

%REQUIRED

%REQUIRED

\section{Introduction}
When the integral equation method is applied to solve a given partial differential equation, one
numerical challenge is the accurate and efficient evaluation of the singular, weakly singular, or hyper
singular integrals representing different potentials in the integral equation reformulations. 
For example, the solutions of a homogeneous elliptic equation (e.g., Laplace, Helmholtz, or Yukawa equations) 
with different types of boundary conditions are often re-expressed as combinations of the single 
layer and double layer potentials with density functions $\rho(z)$ and $\mu(z)$
\begin{equation}\label{eq:SingleDoubleLayer}
\begin{array}{rl}
S \rho (w)= & \int\limits_{\Gamma}^{}G (w,z) \rho (z)\mathrm{d} z, \\
 D \mu (w)= & \int\limits_{\Gamma}^{}\frac{\partial G}{\partial
               \mathbf{n}_{z}} (w,z) \mu (z) \mathrm{d} z,\\
\end{array}
\end{equation}
where $G$ is the free-space Green's function for the underlying elliptic PDE, $z$ is the source point located 
on the boundary $\Gamma$, $w$ is any target point located in the computational domain, and $\mathbf{n}_{z}$  
is the outward normal vector at $z \in \Gamma$. The Green's function $G (w,z)$ is usually a smooth
function when $w$ is away from $z$ on the boundary $\Gamma$, but becomes singular when 
$w \to z$. Therefore, different numerical strategies have to be designed for cases when
$w$ is far away from the boundary, exactly on the boundary, and close to the boundary.

The research topics of developing different numerical integration schemes for evaluating 
the layer potentials at a particular point $w$ for different cases have been extensively
studied. When $w$ is far away from the boundary, classical 
Newton-Cotes or Gaussian quadratures for a general smooth integrand can be applied; 
when $w$ is located on the boundary, special quadrature rules can be designed, 
for example, the trapezoidal rule with end-point corrections 
in \cite{aguilar2002high,alpert1999hybrid,kapur1997high,marin2014corrected,strain1995locally} 
or the generalized Gauss quadrature rules 
in \cite{bremer2010nonlinear,bremer2010universal,yarvin1998generalized};
and when $w$ is close to the boundary, existing techniques include the change of variables to 
remove the principal singularity and the regularized kernel and corrections using asymptotic analysis
\cite{beale2001method,bruno2001fast,duffy1982quadrature,hackbusch1994numerical}. 
We particularly mention the pioneering work in \cite{helsing2008evaluation} which applies 
the Barycentric Lagrange polynomial interpolation formula to derive a globally compensated 
spectrally accurate quadrature rule for evaluation at points close to the boundary 
(also see \cite{barnett2015spectrally}); and the pioneering ``quadrature by expansion" 
(QBX) scheme in \cite{klockner2013quadrature} which derives a partial wave (harmonics) 
expansion valid in a region close to (or even containing points on) the boundary. 
The QBX scheme has been combined with the fast multipole method (FMM) 
in \cite{rachh2017fast} for solving the integral equation reformulation of PDEs.

\begin{figure}[htbp]
\begin{center}
    \includegraphics[width=0.46\linewidth]{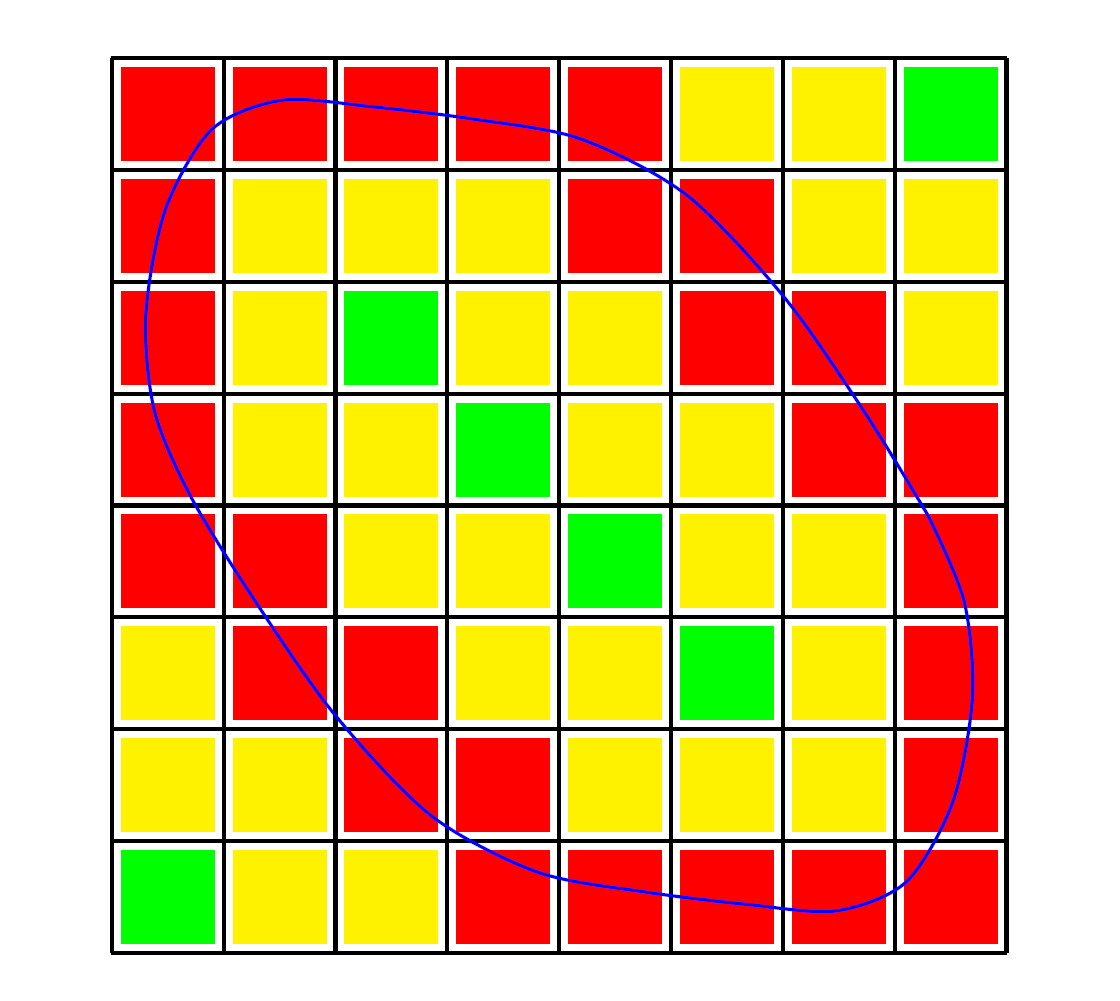}
 \caption[]
  {Different expansions for the leaf boxes in a uniform FMM hierarchical tree structure. 
  {\color{green} Green}: complex polynomial expansion; {\color{yellow} Yellow}: one QB2X for
  the leaf node; {\color{red} Red}: two QB2X required, one for the interior and one for the exterior. }
  \label{fig:fmmgraph}
\end{center}
\end{figure}

In this paper, we introduce new representations of the Laplace layer potentials that are valid 
in the entire leaf (childless) box in the FMM hierarchical tree structure. As the representation 
can be evaluated at any point in the box, the numerical scheme also belongs to the 
class of ``quadrature by expansion" (QBX) schemes to evaluate the 
layer potential integrals. In Fig.~\ref{fig:fmmgraph}, the leaf boxes in a uniform FMM tree 
with $4$ levels are categorized into three groups: The {\color{green} green} box is well 
separated from the boundary source points, from established fast multipole method (FMM) theory 
\cite{greengard1988rapid,greengard1987fast}, the layer potential at each target point 
in the box can be represented by a complex Taylor polynomial expansion which is referred 
to as the local expansion of the green box in the FMM algorithm. Both the {\color{red} red} and 
{\color{yellow} yellow} boxes are not well separated from the boundary, and each of 
the red boxes contains target points located both inside and outside the boundary, 
hence two separate solution representations become necessary, one for the interior 
and one for the exterior. For the red and yellow boxes, contributions from the 
well-separated curved line segments of the layer potential can still be represented 
using the complex local Taylor polynomial expansion, which can be efficiently computed 
using the FMM through the upward and downward passes for the ``far-field" contributions 
of the sources. The numerical difficulty is an accurate and efficient representation 
of the near-field source contributions. 

Using both the complex local Taylor polynomials and plane wave (exponential) functions in 
the basis, we propose a new representation of the 2$D$ layer potentials for the red and 
yellow boxes due to the near-field layer potential source contributions. Combining both the 
far-field and near-field (local) density contributions, the main contribution of this paper 
is that the Laplace layer potentials inside each 2$D$ leaf box of the FMM hierarchical tree 
structure can be represented as the sum of two expansions
\begin{equation}
\label{eq:combined}
	\Re \left( \sum_{k=0}^K c_k (w-w_0)^k + 
	\sum_{p=1}^P \omega_p \frac{e^{\lambda_p \tilde{w}}}{1+\mathrm{i} s'(\tilde{w})} \right).
\end{equation}
Here,  $w= x+ \mathrm{i} y$ is a target point in the leaf box centered at $w_0$, 
the boundary is described by $z=\tx+\mathrm{i} s(\tx)$, $\tilde{w}$ is a point close to $w$ determined by solving the equation 
$z+\mathrm{i} s(z)-w=0$ when $w$ is inside either the red or yellow boxes 
(assuming $s(0)=s'(0)=0$ after proper translations and rotations), the operator $\Re$ takes the real part 
of a complex number, $c_k$ and $\omega_p$ are the complex coefficients of the polynomial and 
plane wave expansions, respectively, the complex number $\lambda_p$ is referred to as 
the node for the exponential (plane wave) expansion, and $K$ and $P$ are respectively the 
numbers of terms in the local Taylor polynomial and plane wave expansions ($P=0$ for the green boxes).

The local complex polynomial expansion only collects the far-field 
contributions and its number of expansion terms $K$ can be analyzed using 
the established error analysis from the classical fast multipole method. 
The plane wave type expansion is derived by applying the Fourier extension
technique to the density and boundary geometry functions followed by 
analytically utilizing the Residue Theorem for complex contour integrals, 
and the number of terms $P$ is the same as the number of terms required 
in the Fourier extensions.
Note that two different types of basis functions are used in the representation of the layer
potential, hence we refer to our approach as the quadrature by two expansions (QB2X). In harmonic 
analysis, the redundant basis functions form a frame \cite{casazza2012finite,duffin1952class}. 
Compared with classical QBX, the QB2X 
representation of the layer potential is valid in a much larger region and allows easier 
analysis of the error and its dependency on the numbers of expansion terms. Another nice feature
of the new representation is that the nonlinear impact of the boundary on the layer potential
becomes explicit in Eq.~(\ref{eq:combined}), providing an analytical tool useful for
other applications.

We organize this paper as follows. In Sec.~\ref{sec:prelim}, we review the classical QBX and 
the well-established Fourier extension technique which form the foundation of the QB2X technique. 
In Sec.~\ref{sec:qbxs}, we derive the new representations for the single and double layer
potentials in the red and yellow boxes using both the complex Taylor polynomial and plane 
wave basis functions. In Sec.~\ref{sec:num}, we present preliminary numerical experiments 
to demonstrate the accuracy of the new representations in the leaf box and to validate our 
analytical results. Finally in Sec.~\ref{sec:summary}, we summarize 
our results and discuss our current work to generalize QB2X to layer potentials
for other types of equations in both two and three dimensions.

\section{Preliminaries}
\label{sec:prelim}
The new quadrature by two expansions (QB2X) technique uses two different basis functions, 
the complex polynomial expansion (also referred to as the local expansion in the fast 
multipole method) and the plane wave expansion using exponential functions. In this section, we 
present (a) the original QBX \cite{klockner2013quadrature} which introduces the complex 
polynomial expansion (or partial waves for the Helmholtz and Yukawa equations) to evaluate 
layer potentials and (b) the Fourier extension technique 
\cite{boyd2002comparison,bruno2007accurate,huybrechs2010fourier} which will provide 
explicit formulas for the plane wave expansion in the new QB2X technique.

\subsection{Quadrature by Expansion: Evaluating Layer Potential Using Complex 
Polynomial Expansion} 
\begin{figure}[htbp]
  \centering
  \subfigure[Analytical layer potential]{
    \includegraphics[width=0.3\linewidth]{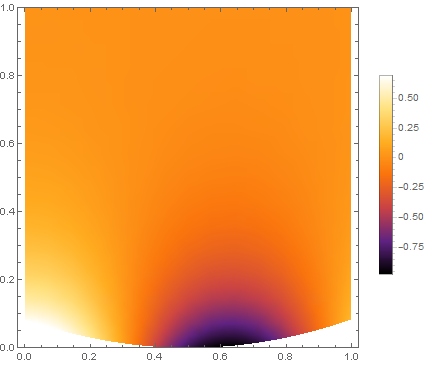}
  }
 \subfigure[$5$ terms expansion ]{
    \includegraphics[width=0.3\linewidth]{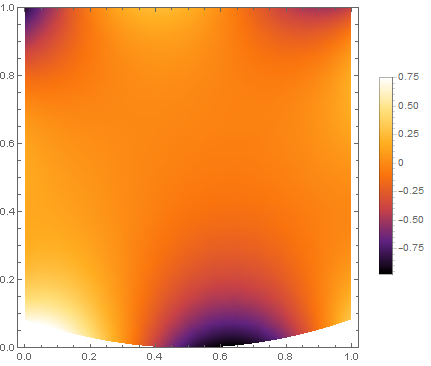}
  }
  \subfigure[$15$ terms expansion]{
    \includegraphics[width=0.3\linewidth]{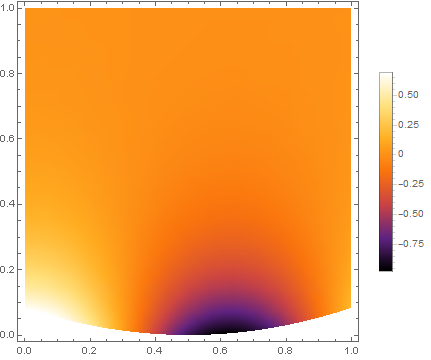}
  }
 \subfigure[$5$ terms error]{
    \includegraphics[width=0.3\linewidth]{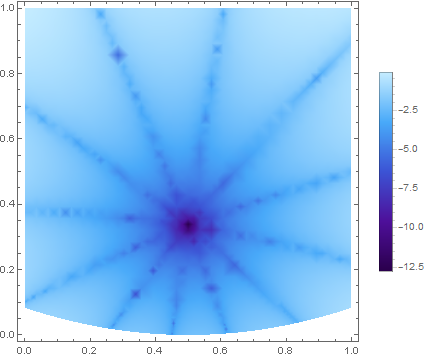}
  }
 \subfigure[$15$ terms error]{
    \includegraphics[width=0.3\linewidth]{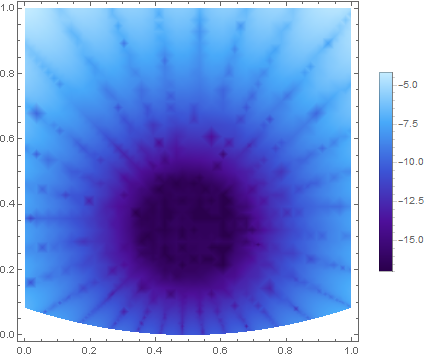}
  }
 \subfigure[$25$ terms error]{
    \includegraphics[width=0.3\linewidth]{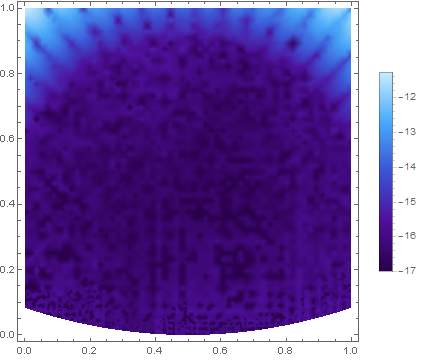}
  }
  \hfill
  \caption[]
  {An implementation of QBX.  (a): Analytical layer potential. 
   (b) and (c): approximation using $K=5$ and $K=15$ terms, respectively.
	(d), (e), and (f): $\log_{10}$ errors for $K=5$, $K=15$, and $K=25$. 
 }
  \label{fig:qbxgraph}
\end{figure}

Assuming both the density function and boundary curve are sufficiently smooth, to evaluate
the singular, near-singular, or hyper-singular layer potentials, in \cite{klockner2013quadrature},
it was observed that the layer potentials are smooth functions on either side of the
boundary, and the integrand singularity is only associated with the non-smoothness
across the boundary. Therefore the polynomial expansion of the Laplace layer potential
centered at a point either in the interior or exterior of the boundary is valid at least locally.
In Fig.~\ref{fig:qbxgraph}, we consider a Laplace layer potential explicitly given by
$\Re \left( e^{ \mathrm{i} 5 z} \right) = \Re \left( e^{\mathrm{i} 5
    (x+ \mathrm{i}y)} \right) $ in a box $0<x,y<1$ where $z=x+ \mathrm{i}y$ is the complex 
variable. The box contains part of the boundary given by the equation $\tilde{y}=s(\tx)=\frac13 (\tx-\frac12)^2$.
We assume the polynomial expansion is in the form $ \Re \left( \sum_{k=0}^K c_k (w-w_0)^k \right) $ 
centered at $w_0=\frac12+\frac13 \mathrm{i}$. We neglect the numerical errors when evaluating 
the expansion coefficients, i.e., the coefficients are derived exactly. 
In (a), we plot the analytical layer potential. In (b) and (c), we present the polynomial 
approximation of the layer potential using $K=5$ and $K=15$ terms in the expansion, and 
in (d), (e), and (f), we plot the $\log_{10}$ errors of the representation when 
$K=5$, $K=15$, and $K=25$, respectively. Clearly, when the number of expansion terms 
$K$ increases, the error decreases and the representation becomes valid in a much 
larger region.

Using standard error analysis from the FMM theory, assuming the coefficients are 
computed accurately for the well-separated green boxes in Fig.~\ref{fig:fmmgraph}, 
the local expansion in the form $ \Re \left( \sum_{k=0}^K c_k (w-w_0)^k \right) $ 
can achieve $6$-digits accuracy when $K$ is approximately $18$ and $12$-digits 
accuracy when $K=36$. The task of evaluating the layer potentials at 
points in the red and yellow boxes becomes more complicated.  
In \cite{klockner2013quadrature}, the choice of the local complex polynomial 
expansion center, the degree of the polynomial, and the quadrature schemes 
for computing the expansion coefficients are numerically studied to provide 
guidelines on these parameter selections. In \cite{epstein2013convergence}, 
estimates for the rate of convergence of these local expansions are derived,
which can be used to analyze the approximation error of the local Taylor
expansions in a leaf box. In existing FMM+QBX implementation \cite{rachh2017fast}, 
for evaluation points in one leaf node, several complex 
polynomial expansions may have to be formed with different expansion centers 
and degrees, unless the FMM hierarchical tree oversamples the density on the boundary 
or the solution. We also mention that in existing QBX implementations, the curvature of
the boundary does not explicitly appear in the formulas. This will be addressed
in the new QB2X in Sec.~\ref{sec:qbxs}.

\subsection{Fourier Extension: Approximation Using Exponentials}
Compared with a polynomial basis, the exponential expansions, if it can be derived
accurately and efficiently, may show better numerical properties in efficiency.
One example is the translations in the FMM algorithms. When the exponential
(plane wave) expansions are used, the translations become diagonal and
the number of operations is reduced from the polynomial expansion's $O(K^2)$ 
to the plane wave expansion's $O(K)$ when $K$ terms are used in both 
expansions \cite{cheng2006wideband,crutchfield2006remarks,greengard1997new}. 
Unfortunately, deriving the optimal exponential expansion for a general function requires
nonlinear optimization, and the uniqueness of the  solution is not 
guaranteed. However, in some particular cases, a good exponential expansion
approximation can be derived. For example, when the function is smooth and 
periodic, then the Fourier series can be computed efficiently and the expansion 
converges rapidly. Another example is when the inverse integral transform 
of the function is available, then the exponential expansion problem becomes 
an integration problem, and the weights and nodes of the exponential expansions
can be computed using the generalized Gauss quadrature method 
\cite{ma1996generalized,yarvin1998generalized}.  
\begin{figure}[htbp]
  \centering
  \subfigure[$f (x)=x$]{
    \includegraphics[width=0.46\linewidth]{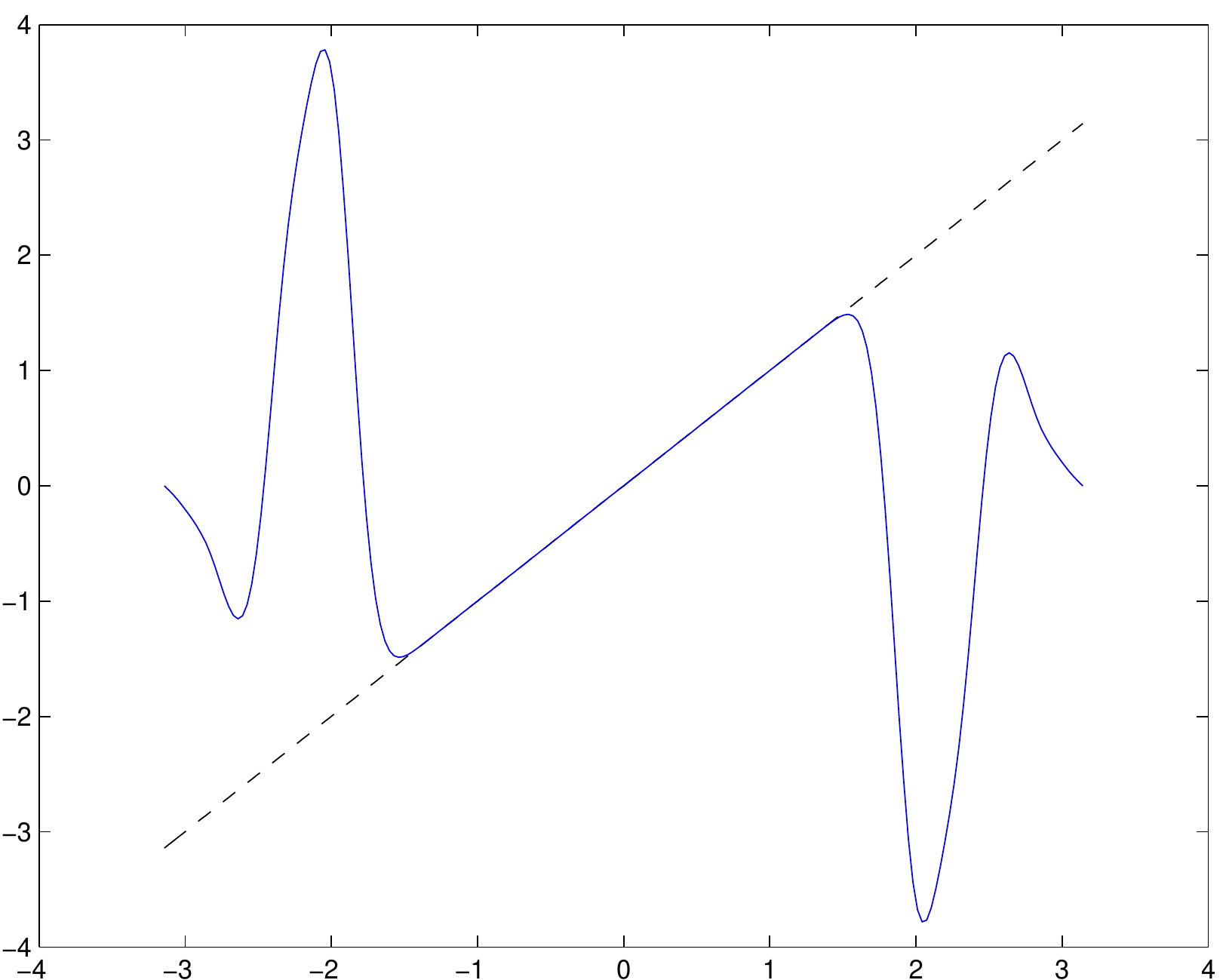}
  }
 \subfigure[$f (x)=x$]{
    \includegraphics[width=0.46\linewidth]{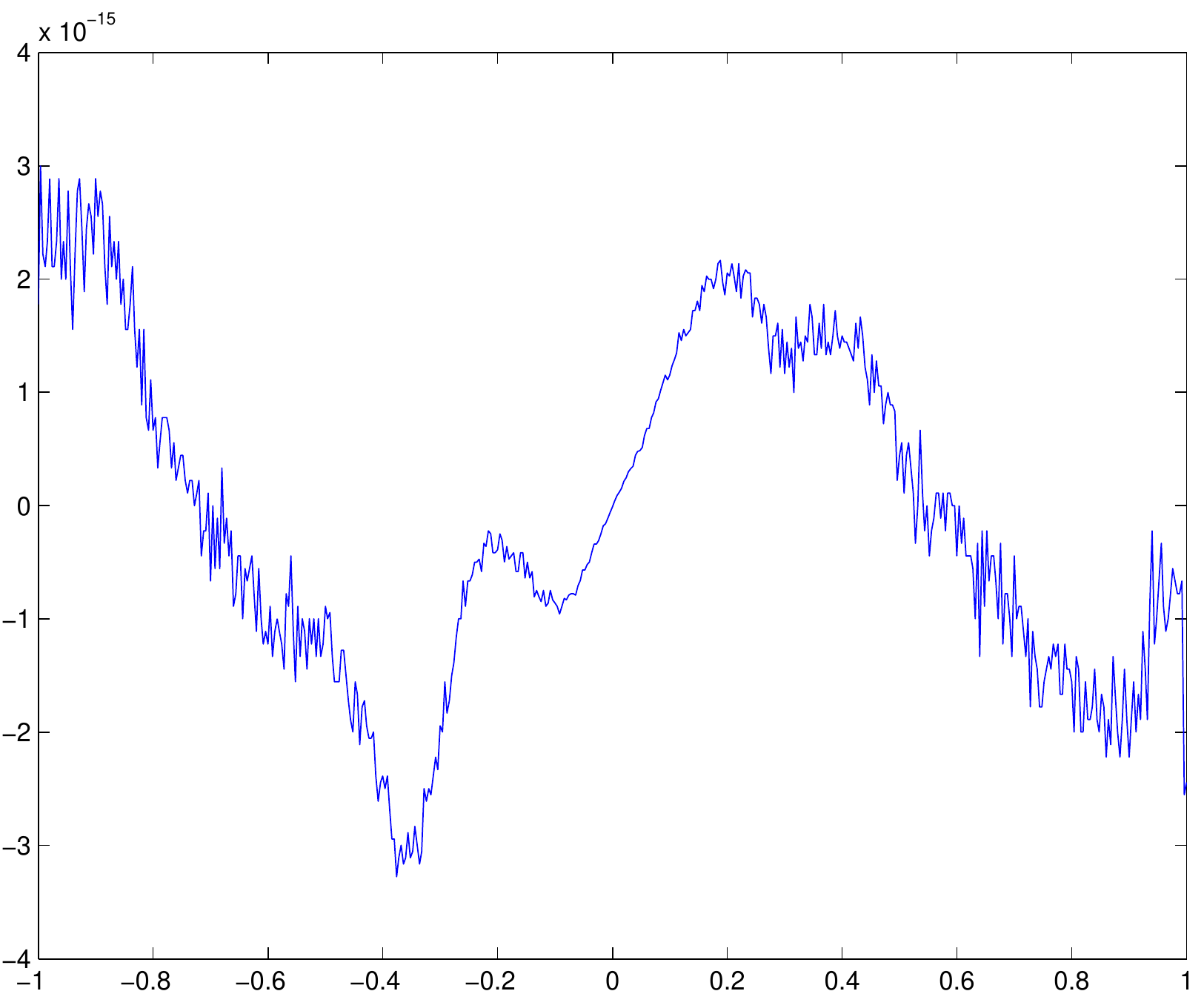}
  }
  \subfigure[$f (x)=2x^2-1$]{
    \includegraphics[width=0.46\linewidth]{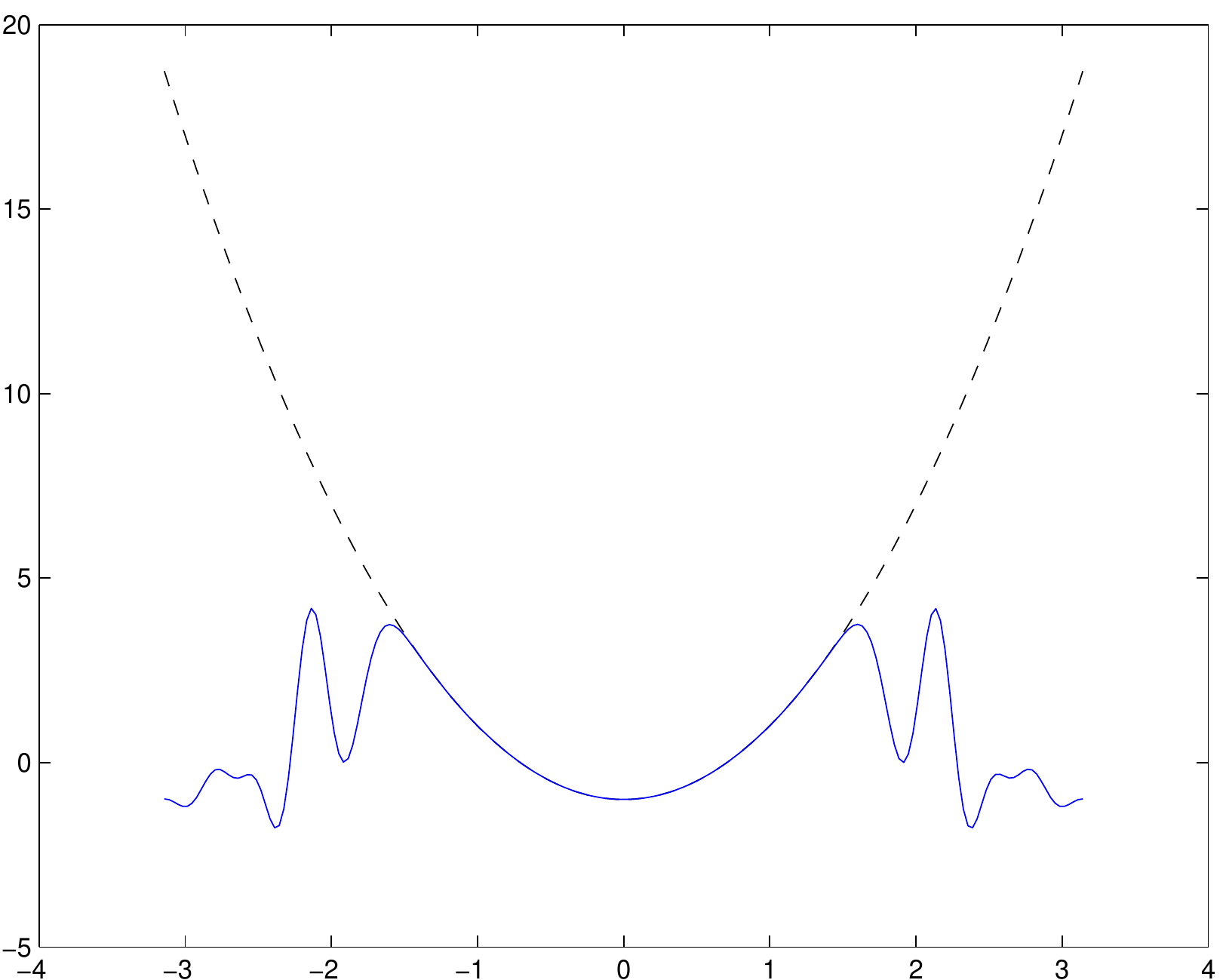}
  }
 \subfigure[$f (x)=2x^2-1$]{
    \includegraphics[width=0.46\linewidth]{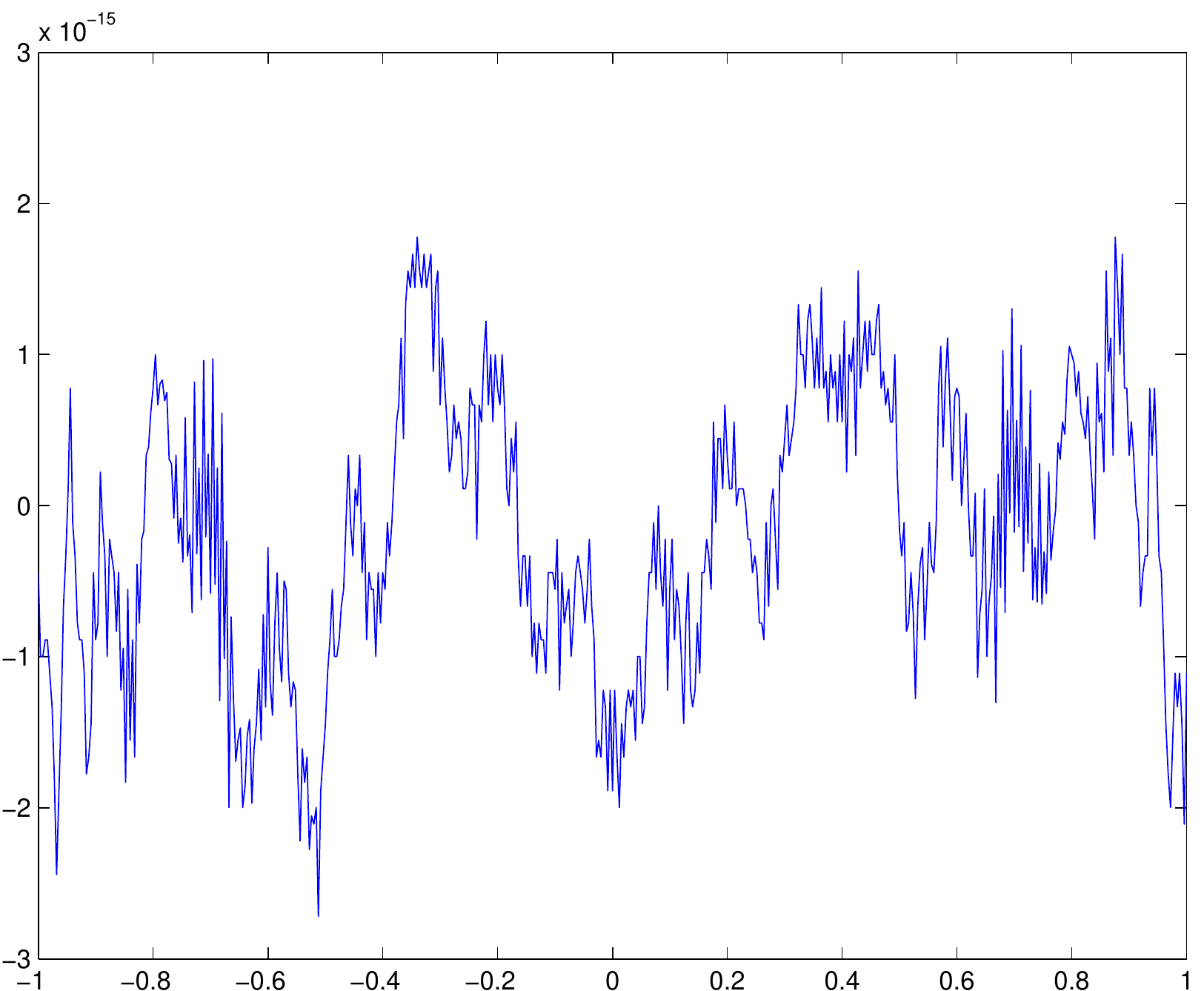}
  }
  \hfill
  \caption[]
	{(a) and (c): Fourier extension $g(x)=\sum_{p=-P}^P c_p e^{\mathrm{i} p x}$ (solid line) of the 
	given function $f$ (dashed line), $P=30$. (b) and (d): approximation errors on $[-1,1]$. 
 }
  \label{fig:Fourierextension}
\end{figure}
Using the Fourier series to approximate a non-periodic function is also
a well-studied topic. Consider a non-periodic function $f$ on $[-1,1]$, 
applying the Fourier extension technique, a suitable periodic function $g$ on 
a larger domain $[-T,T]$ is computed stably, so the Fourier series expansion of 
$g$ matches $f$ on the interval $[-1,1]$, 
see \cite{boyd2002comparison,bruno2007accurate,huybrechs2010fourier} and
references therein. We have implemented the scheme in \cite{huybrechs2010fourier},
which solves the least square optimization problem to compute an accurate Fourier 
series representation of a smooth function defined on $[-1,1]$. 
In Fig.~\ref{fig:Fourierextension}, we present the computed Fourier series with
fundamental period $2 \pi$ for the two Chebyshev basis polynomials $T_1(x)=x$ and
$T_2(x)=2 x^2 -1$.  The Fourier series approximation on $[-1,1]$ achieves 
machine precision accuracy for both cases. Using the Fourier extension technique, a translation 
matrix can be precomputed to map the commonly used Chebyshev or other orthogonal polynomial 
basis functions to the Fourier basis functions that are periodic in a larger domain. 

\section{Quadrature by Two Expansions: Combining Complex Polynomial and Plane Wave Expansions}
\label{sec:qbxs}
\subsection{Problem Setup}
\begin{figure}[htbp]
\begin{center}
	\includegraphics[width=0.66\linewidth]{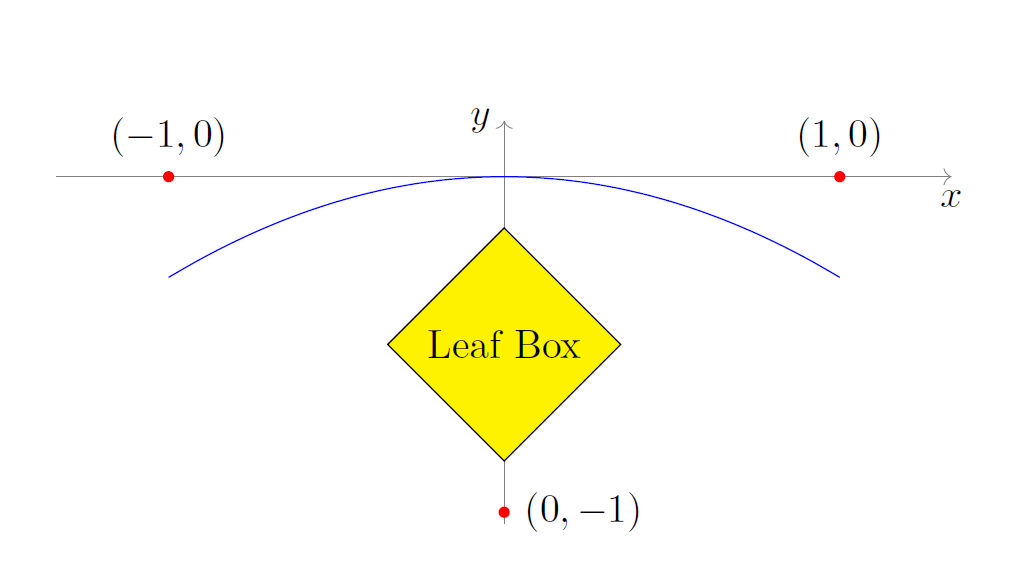}
 \caption[]
  {A leaf box close to the boundary. }
  \label{fig:setting}
\end{center}
\end{figure}

As all the far-field layer potential density contributions can be accurately and 
efficiently computed using the fast multipole method, in this section, we focus
on the near-field (local) density contributions to the red and yellow boxes in
Fig.~\ref{fig:fmmgraph}. We assume the boundary is parametrically described by
$z=\tx +\mathrm{i} s(\tx)$,$-1\leq \tx \leq 1$ and $s (0)=0, s'(0)=0$ 
after proper scaling, translation, and rotation as in Fig.~\ref{fig:setting}. 
One of the yellow leaf box is shown and we assume the two end points $(-1,s(-1))$ 
and $(1, s(1))$ are well-separated from this leaf box with center $w_0$.
Following previous work in \cite{rachh2017fast}, we assume both the density function 
$\rho(\tilde{x})$ and boundary $s(\tilde{x})$ are well resolved by polynomials for 
$-1 \leq \tilde{x} \leq 1$.
Note that the degree of the polynomials to approximate $s(\tilde{x})$ and $\rho(\tilde{x})$
can be different. Next, we present the detailed representation formulas for the single and 
double layer potentials 
\begin{equation}
\begin{array}{rl}
S \rho (w)= & \int_{-1}^{1} G (w,z) \rho (\tilde{x}) |z'(\tilde{x})| \mathrm{d} \tilde{x}, \\
D \rho (w)= & \int_{-1}^{1} \frac{\partial G}{\partial \mathbf{n}_{z}} (w,z) \rho (\tilde{x}) |z'(\tilde{x})| \mathrm{d} \tilde{x} \\
\end{array}
\end{equation}
where $w= x+\mathrm{i} y$ is the target point in the yellow or red leaf box, 
$z=\tilde{x}+\mathrm{i} s(\tilde{x})$ is the source on the boundary, and 
$|z'(\tilde{x})|=\sqrt{1+ s' (\tilde{x})^2}$.

\subsection{Laplace Double Layer Potential}
We start from the double layer potential to avoid the branch cut analysis of the 
$\log$ function in the kernel. To evaluate the double layer potential at the target point 
$w= x + \mathrm{i} y$ in the yellow or red leaf box, the contribution from the source density function 
$\rho(\tilde{x})$ defined on the boundary segment $z =\tilde{x} + \mathrm{i} s (\tilde{x})$,  
$-1 \leq \tilde{x} \leq 1,$ becomes
\begin{equation}
\begin{array}{rl}
DLP (w=x+ i y)= &  \int\limits_{-1}^1  \frac{\partial G (w,z)}{\partial
             \mathbf{n}_z} \rho (\tx) \left| z' (\tilde{x}) \right| \mathrm{d} \tilde{x}\\
 = &  \frac{1}{2 \pi} \int\limits_{-1}^1   \frac{ \left\langle x-\tilde{x}, y-s
     (\tilde{x}) \right\rangle}{ (x-\tilde{x})^2+ (y-s (\tilde{x}))^2}
     \cdot \frac{\left\langle s '(\tilde{x}),-1 \right\rangle}{\sqrt{1+
     s' (\tilde{x})^2}} \rho (\tilde{x}) \sqrt{1+ s' (\tilde{x})^2}
     \mathrm{d} \tilde{x}\\
=&  \frac{1}{2\pi} \int\limits_{-1}^1 \frac{(x-\tilde{x}) s'
   (\tilde{x})- (y- s (\tilde{x}))}{(x-\tilde{x})^2+ (y- s
   (\tilde{x}))^2} \rho (\tilde{x}) \mathrm{d} \tilde{x}
\end{array}    
\end{equation}
where $\frac{\partial G}{ \partial \mathbf{n}_z}= \frac{1}{2\pi}
\frac{w-z}{\lVert w-z\rVert } \cdot \mathbf{n}_z (\tilde{x})$,
$\mathbf{n}_z (\tilde{x})= \frac{\left\langle s '(\tilde{x}),-1 \right\rangle}{\sqrt{1+s' (\tilde{x})^2}} $. 
The curvature of the boundary at $z=(\tilde{x},s (\tilde{x}))$ is given by $\kappa= \frac{s''
  (\tilde{x})}{(1+\left| f' \right|^2)^{\frac{3}{2}}}$.
In order to apply the complex contour integral theory and Residue Theorem, using 
$$\frac{1}{w-z}= \frac{\overline{w-z}}{(w-z) \overline{(w-z)}}=
\frac{(x-\tilde{x})- \mathrm{i} (y- s (\tx))}{ (x-\tilde{x})^2+ (y- s(\tilde{x}))^2},$$ 
we get $$\frac{(x-\tilde{x}) s'(\tilde{x})}{(x-\tilde{x})^{2}+ (y- s(\tilde{x}))^2} 
= \Re \left(\frac{1}{w-z} s' (\tilde{x}) \right), \quad
\frac{- (y- s (\tilde{x}))}{(x-\tilde{x})^{2}+ (y- s(\tilde{x}))^2} =\Im \left(\frac{1}{w-z} \right).$$ 
The double layer potential becomes 
\begin{equation}
\begin{array}{rl}
DLP (w)= &  \frac{1}{2\pi} \int\limits_{-1}^1 \Re (\frac{1}{w-z})s'
           (\tilde{x}) \rho (\tilde{x}) \mathrm{d} \tilde{x}+
           \frac{1}{2\pi} \int\limits_{-1}^1  \Im (\frac{1}{w-z})
           \rho (\tilde{x}) \mathrm{d} \tilde{x}\\
=&  \frac{1}{2\pi} \Re\int\limits_{-1}^1  \frac{s'
           (\tilde{x}) \rho (\tilde{x})}{w- (\tilde{x}+
   \mathrm{i} s (\tilde{x}))} \mathrm{d} \tilde{x}+
           \frac{1}{2\pi} \Im\int\limits_{-1}^1  \frac{ \rho (\tilde{x})}{w-
(\tilde{x}+\mathrm{i} s (\tilde{x}))} \mathrm{d} \tilde{x}.  \\
\end{array}
\end{equation}
As $s(\tilde{x})$ and $\rho(\tilde{x})$ are resolved by polynomials, therefore
both terms are in the form of the complex integral
\begin{equation}
\int\limits_{-1}^1  \frac{f(\tilde{x})}{(\tilde{x}+ \mathrm{i} s (\tilde{x}))-w} \mathrm{d} \tilde{x}
\label{eq:dlp}
\end{equation}
where $f(\tilde{x})$ is a polynomial defined for $\tilde{x} \in [-1,1]$.

\subsubsection{$s(\tilde{x})=0$}
We first consider the case when $s(\tilde{x})=0$ to simplify the discussions and formulas, i.e.,
$z$ is on the straight line segment connecting $(-1,0)$ and $(1,0)$. In this case, Eq.~(\ref{eq:dlp}) 
becomes $$\int\limits_{-1}^1 \frac{f (z)}{z-w} \mathrm{d} z.$$ 
As $f(z)$ is a function defined on the real line segment, the Fourier extension technique can be 
applied using the precomputed translation operator from the polynomial basis to the 
Fourier series, so $$f (z) \approx \sum\limits_{p=-P}^{P} \omega_p e^{\mathrm{i} p z}$$ 
and Eq.~(\ref{eq:dlp}) becomes 
\begin{equation}
\begin{array}{rl}
\int\limits_{-1}^1 \frac{f (z)}{z-w}\mathrm{d} z & \approx
 \int\limits_{-1}^{1} \frac{1}{z-w} ( \sum\limits_{p=-P}^{P} \omega_p e^{\mathrm{i} p z}) \mathrm{d} z \\
 & =\int\limits_{-1}^{1}  \sum\limits_{p=0}^{P}\frac{1}{z-w} \omega_pe^{\mathrm{i} p z} \mathrm{d} z+  
  \int\limits_{-1}^{1}  \sum\limits_{p=-P}^{-1}\frac{1}{z-w} \omega_pe^{\mathrm{i} p z} \mathrm{d}z\\
& = I_1+I_2,\\
\end{array}
\label{eq:I1I2}
\end{equation}
where $I_1$ only contains the non-negative $p$ frequencies and $I_2$ contains the negative ones.

\begin{figure}[htbp]
\begin{center}
\includegraphics[width=0.86\linewidth]{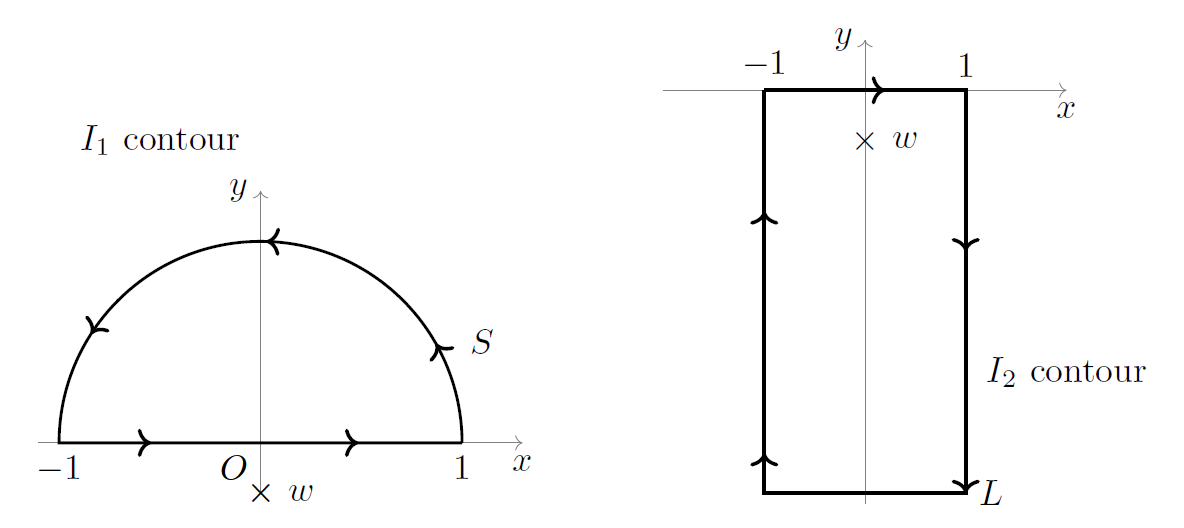}
 \caption[]{The upper (left) and lower (right) contours for $I_1$ and $I_2$, respectively. }
  \label{fig:contours}
\end{center}
\end{figure}
We first study $I_1$ as part of the contour integral 
$$ \sum\limits_{p=0}^{P} \int\limits_{C}^{}  \frac{1}{z-w} \omega_pe^{\mathrm{i} p z}\mathrm{d}z $$
where the contour $C$ is shown in the left plot of Fig.~\ref{fig:contours}. 
It consists of the line segment from $-1$ to $1$ and the semi-circle (denoted by $S$)  
on the upper half plane.
As the integrand is analytic inside the contour and the pole is located at $w$ outside 
the contour, by the Residue Theorem, we have 
$$ \sum\limits_{p=0}^{P} \int\limits_{C}^{}  \frac{1}{z-w} \omega_pe^{\mathrm{i} p z}\mathrm{d}z=0. $$
Consequently, 
\begin{equation}
I_1= \int\limits_{-1}^{1} \sum\limits_{p=0}^{P}\frac{1}{z-w} \omega_pe^{\mathrm{i} p z} 
\mathrm{d}z = -\sum\limits_{p=0}^{P} \int_{S} \frac{1}{z-w} \omega_pe^{\mathrm{i} p z}\mathrm{d}z. 
\end{equation}
Unlike the line segment from $-1$ to $1$, the semi-circle $S$ is well-separated from the leaf box containing
$w$ and the contribution from the density defined on the semi-circle can be collected into
a local expansion as in
\begin{equation}
\begin{array}{rl}
 I_1= &-\sum\limits_{p=0}^{P} \int\limits_{S}^{}  \frac{1}{(z-w_{0})-(w-w_0)} \omega_pe^{\mathrm{i} p z} \mathrm{d}z
 = -\sum\limits_{p=0}^{P} \int\limits_{S}^{}
     \frac{1}{(z-w_{0})(1-\frac{w-w_0}{z-w_0})} \omega_p e^{\mathrm{i}
     p z} \mathrm{d} z \\
 \approx &-\sum\limits_{p=0}^{P} \int\limits_{S}^{}  \frac{1}{(z-w_{0})}\sum\limits_{k=0}^{K} (\frac{w-w_0}{z-w_0})^k  \omega_p 
  e^{\mathrm{i} p z} \mathrm{d}z = -\sum\limits_{k=0}^{K}  c_k (w-w_0)^k \\
\end{array}
\end{equation}
where $w_{0}$ is the center of the leaf box, $c_k= \left( \sum\limits_{p=0}^{P}   \omega_p \int\limits_{S}^{}
          \frac{ e^{\mathrm{i} p z}}{(z-w_0)^{k+1}} \mathrm{d}z \right) $ are the local expansion
coefficients and the number of terms $K$ is controlled by the decay rate of $|(\frac{w-w_0}{z-w_0})^k|$ 
which can be easily estimated for the given contour and target leaf box following the 
standard fast multipole method analysis. Introducing $r_{max}=max_{ \{z,w \} } |\frac{w-w_0}{z-w_0}|$ and 
noting that when $p$ is non-negative, $|e^{i p z}|$ decreases when $z$ on the contour moves away from 
the line segment, a very loose estimate of the truncation error is given by 
$c \left( \sum_{p=0}^{P} | \omega_p|  \right) r_{max}^{K+1}$ for some constant $c$. When $s(\tilde{x})=0$, 
it is straightforward to verify that $K=9$, $18$, $27$, and $36$ will provide results with 
at least $3$, $6$, $9$, and $12$ accurate digits, respectively.
Also, the numerical stability issues associated with exponentially growing $|e^{i p z}|$ values 
can be avoided and the local expansion coefficients can be computed accurately using standard
quadrature rules. 
Therefore, the non-negative modes can be represented as a local complex Taylor polynomial expansion. 
Unfortunately, the negative $p$ frequencies cannot be computed using this contour, as the 
function $e^{i p z}$ grows exponentially when $z$ on the contour moves away from the real axis 
for $p < 0$.

To compute $I_2$, a different contour $C$ on the lower half complex plane has to be chosen, 
a sample contour is show on the right plot of Fig.~\ref{fig:contours}. It consists of a rectangle
with one side being the line segment from $-1$ to $1$, and a sufficiently large (can be $\infty$) 
constant $L$ is introduced to determine the length of the other side. We denote the part of 
the contour consisting of the three other sides of the rectangle by $S$. 
Note that $e^{i p z}$ decays exponentially when $z$ on $S$ moves away from the real axis for $p<0$. 
As $w$ is inside the contour, applying the Residue Theorem, for each negative frequency $p$, 
we have 
$$ \int\limits_{-1}^{1} \frac{1}{z-w} \omega_pe^{\mathrm{i} p z} \mathrm{d}z 
  +  \int\limits_S^{} \frac{1}{z-w} \omega_pe^{\mathrm{i} p z} \mathrm{d}z 
 = - 2\pi \mathrm{i} Res[\frac{1}{z-w} \omega_pe^{\mathrm{i} p z} ,w ].
$$
Therefore, we can compute $I_2$ using
\begin{equation}
\begin{array}{rl}
I_2= &\int\limits_{-1}^{1}\sum\limits_{p=-P}^{-1}\frac{1}{z-w} \omega_pe^{\mathrm{i} p z} \mathrm{d}z  \\
 = &-2\pi \mathrm{i} Res[\sum\limits_{p=-P}^{-1}\frac{1}{z-w} \omega_pe^{\mathrm{i} p z} ,w ] - \int\limits_S^{} \sum\limits_{p=-P}^{-1}\frac{1}{z-w} \omega_pe^{\mathrm{i} p z} \mathrm{d}z \\
 = &-2\pi\mathrm{i} \sum\limits_{p=-P}^{-1} \omega_pe^{\mathrm{i} p w} - \int\limits_S^{} \sum\limits_{p=-P}^{-1}\frac{1}{z-w} \omega_pe^{\mathrm{i} p z} \mathrm{d}z \\
\approx &-2\pi \mathrm{i} \sum\limits_{p=-P}^{-1} \omega_pe^{\mathrm{i} p w}
    -\sum\limits_{k=0}^{K} c_k (w-w_0)^k \\
\end{array}
\end{equation}
where the local expansion (second summation in the formula) coefficients are given by
$$c_k= \sum\limits_{p=-P}^{-1}  \omega_p \int\limits_{S}^{} 
\frac{ e^{\mathrm{i} p z}}{(z-w_0)^{k+1}} \mathrm{d}z$$
which are derived using the same separation of variables as in $I_1$. 
As the leaf box is well-separated from $S$, the number $K$ of the local 
polynomial expansion can be determined using the same FMM error analysis as
in $I_1$ and we skip the details.

Combining $I_1$, $I_2$, and the far-field density contributions for this special case, 
we conclude that the double layer potential in the leaf box can be represented as a 
combination of the local complex Taylor polynomial expansion and plane wave expansion as in 
Eq.~(\ref{eq:combined}). The number of terms in the local polynomial expansion is determined 
by standard FMM error analysis as all the involved contributions are well-separated from
the leaf box. The number of terms in the plane wave expansion is the same as that
in the Fourier extension of the density and boundary functions for $-1<\tilde{x}<1$.

\vspace{0.05in}
{\noindent \bf Comment on $r_{max}$:} Smaller $r_{max}$ values are possible by including a larger
portion of the boundary when computing $I_1$ and $I_2$, e.g., by also including contributions from
the second nearest neighbors \cite{nabors1991fastcap}, at the cost of more terms in the Fourier 
extension. The balance of the numbers of terms $P$ in the exponential expansion and $K$ in the
polynomial expansion is related with the optimal discretization strategies when 
generating the FMM adaptive tree, which is currently being studied.

\vspace{0.05in}
{\noindent \bf Comment on the contours:} We mention that the choice of the contour
is not unique. Other contours can also be used. However changing the contours will not
change the values of the expansion coefficients. It will only change the accuracy and 
efficiency of the numerical integration scheme for computing these values and the 
estimated number $K$ for truncating the local Taylor expansion.
Our numerical experiments in Sec.~\ref{sec:num} show that the current choice 
allows accurate and efficient computations of these coefficients and provides
acceptable bounds for $K$. Finding the ``optimal" contours is still a challenging 
task and is being studied.

\subsubsection{Curved Line: Role of Boundary Geometry}
\label{sec:doubleCurve}
Next we consider the role of the boundary geometry in the representation. We assume $s(\tilde{x})$ 
is a polynomial of $\tilde{x}$, $s(\tilde{x})=s'(\tilde{x})=0$ at $\tilde{x}=0$, and 
$|s (\tilde{x})| \ll |\tilde{x}|$. The meaning of the last assumption is
that the boundary is sufficiently resolved by $s(\tilde{x})$, so that the boundary curve
in each leaf box is reasonably close to flat. In this case, the double layer potential 
representation will depend on the boundary geometry described by the polynomial $s(z)$.
When the degree of the polynomial $s(z)$ is less than $5$, the roots of  
$z+\mathrm{i} s (z)-w=0$ can be derived analytically, and the one close to $w$ will 
be denoted as $\tilde{w}$. The assumption that $|s (z)| \ll |z|$  for $|z|$ small also 
implies that all other roots will be very far away from the leaf box and the region 
enclosed by the contours. Similar to the $s(\tilde{x})=0$ case, we replace $f(z)$ by 
its Fourier series expansion, and separate the integral in Eq.~(\ref{eq:dlp}) into two parts
\begin{equation}
\begin{array}{rl}
 \int\limits_{-1}^1 \frac{f (z)}{(z+\mathrm{i} s (z))-w}\mathrm{d} z = & 
    \int\limits_{-1}^{1} \frac{1}{z+\mathrm{i} s(z)-w} (\sum\limits_{p=-P}^P  \omega_p e^{\mathrm{i} p z}) \mathrm{d} z\\
 = & \int\limits_{-1}^{1} \sum\limits_{p=0}^P \frac{1}{z+\mathrm{i} s(z)-w}  \omega_p e^{\mathrm{i} p z} \mathrm{d} z+  
    \int\limits_{-1}^{1}  \sum\limits_{p=-P}^{-1}\frac{1}{z+\mathrm{i} s(z)-w}  \omega_p e^{\mathrm{i} p z} \mathrm{d}z\\
=& I_1+I_2\\
\end{array}
\end{equation}
 
To evaluate $I_1$, we use the same contour on the left of Fig.~\ref{fig:contours}.
Note that the semi-circle (denoted by $S$) is well-separated from the leaf box,
the integrand is analytic in the region enclosed by the contour, and $e^{\mathrm{i} p z}$
decays exponentially when $p>0$ and $\Re (z) \to + \infty$. Therefore no numerical stability 
issues will appear and $I_1$ can be approximated using the same strategy as in 
the $s(\tilde{x})=0$ case by a complex polynomial expansion 
\begin{equation}
\begin{array}{rl}
 I_1= &-\sum\limits_{p=0}^{P} \int\limits_{S}^{}  \frac{1}{(z+\mathrm{i} s(z)-w_{0})-(w-w_0)}  \omega_p e^{\mathrm{i} p z} \mathrm{d}z \\
 = & -\sum\limits_{p=0}^{P} \int\limits_{S}^{}
     \frac{1}{(z+\mathrm{i} s(z)-w_{0})(1-\frac{w-w_0}{z+\mathrm{i} s(z)-w_0})}  \omega_p e^{\mathrm{i}
     p z} \mathrm{d} z \\
 \approx &-\sum\limits_{p=0}^{P} \int\limits_{S}^{}  \frac{1}{(z+\mathrm{i} s(z)-w_{0})}\sum\limits_{k=0}^{K} 
  (\frac{w-w_0}{z+\mathrm{i} s(z)-w_0})^k  \omega_p e^{\mathrm{i} p z} \mathrm{d}z = -\sum\limits_{k=0}^{K}  c_k (w-w_0)^k \\
\end{array}
\end{equation}
where $w_{0}$ is the center of the leaf box, $$c_k= \left( \sum\limits_{p=0}^{P}   \omega_p \int\limits_{S}^{}
\frac{ e^{\mathrm{i} p z}}{(z+\mathrm{i} s(z)-w_0)^{k+1}} \mathrm{d}z \right) $$ are the local expansion
coefficients, and the number of terms $K$ can be determined by standard FMM error analysis using the
ratio $r_{max}=max_{ \{z,w \} } |\frac{w-w_0}{z+\mathrm{i} s(z)-w_0}|$ for $z$ on the contour and $w$ 
in the leaf target box. Note that when $|s(\tilde{x})| \ll |\tilde{x}|$, the ratio $r_{max}$ only changes
slightly when compared with the $s(\tilde{x})=0$ case.

We also apply the same contour on the right of Fig.~\ref{fig:contours} to evaluate $I_2$.
As the root $\tilde{w}$ is inside the contour so we factor $z+\mathrm{i} s(z)-w = (z-\tilde{w}) g(z)$. 
Simple algebra will show that 
$g(\tilde{w})= 1+ \mathrm{i} s'(\tilde{w})$ and the Residue Theorem becomes
$$ \int\limits_C \frac{1}{z+\mathrm{i} s(z)-w}  \omega_p e^{\mathrm{i} p z} \mathrm{d}z 
  = - 2\pi \mathrm{i} Res[\frac{1}{z+\mathrm{i} s (z)- w}  \omega_p e^{\mathrm{i} p z}, \tilde{w} ] 
  = - 2\pi \mathrm{i}  \omega_p \frac{e^{\mathrm{i} p \tilde{w}}}{1+\mathrm{i} s'(\tilde{w})}.
$$
We can therefore represent $I_2$ using the following formula, 
\begin{equation}
\begin{array}{rl}
I_2= &\int\limits_{-1}^{1}\sum\limits_{p=-P}^{-1}\frac{1}{z+\mathrm{i} s(z)-w}  \omega_p e^{\mathrm{i} p z} \mathrm{d}z  \\
 = &-2\pi \mathrm{i} Res[\sum\limits_{p=-P}^{-1}\frac{1}{z+\mathrm{i} s(z)-w}  \omega_p e^{\mathrm{i} p z},\tilde{w}] - 
  \int\limits_S^{} \sum\limits_{p=-P}^{-1}\frac{1}{z+\mathrm{i} s(z)-w} \omega_pe^{\mathrm{i} p z} \mathrm{d}z \\
 = &-2\pi\mathrm{i} \sum\limits_{p=-P}^{-1}  \omega_p \frac{e^{\mathrm{i} p \tilde{w}}}{1+\mathrm{i} s'(\tilde{w})} 
  - \int\limits_S^{} \sum\limits_{p=-p}^{-1}\frac{1}{z+\mathrm{i} s(z) -w}  \omega_p e^{\mathrm{i} p z} \mathrm{d}z \\
  \approx &-2\pi \mathrm{i} \sum\limits_{p=-P}^{-1} \omega_p \frac{e^{\mathrm{i} p \tilde{w}}}{1+\mathrm{i} s'(\tilde{w})}
    -\sum\limits_{k=0}^{K} c_k (w-w_0)^k \\
\end{array}
\end{equation}
where the local expansion coefficients are given by
$$c_k= \sum\limits_{p=-P}^{-1}  \omega_p \int\limits_{S}^{} 
  \frac{ e^{\mathrm{i} p z}}{(z+\mathrm{i} s(z) -w_0)^{k+1}} \mathrm{d}z,$$
and the number of terms $K$ can be estimated using the ratio $r_{max}$ and standard FMM error analysis. 
This formula shows that the representation depends nonlinearly on the geometry in three different ways,
the additional $s'(\tilde{w})$ term in the denominator of the (slightly modified) plane wave expansion,
the nonlinear dependency when finding the root $\tilde{w}$ using the polynomial equation 
$z+\mathrm{i} s(z)-w =0$, and the ratio $r_{max}$ which depends on the function $s(z)$. 
When the degree of the polynomial $s(z)$ is less than $5$, an analytical
formula is available to express $\tilde{w}$ explicitly as a function of $w$. When the degree is higher,
an asymptotic expansion can be derived to approximate $\tilde{w}$ using $w$ when assuming 
$|s (\tilde{x})| \ll |\tilde{x}|$. The resulting mathematical formulas in the QB2X therefore reveal 
the solution dependency on the geometry and can be useful tools in PDE analysis.

\vspace{0.05in}
{\noindent \bf Comment on $|s (\tilde{x})| \ll |\tilde{x}|$:} It is worth mentioning that even without 
this assumption, most parts in the analysis are still valid. It is therefore possible to apply
the combined (slightly modified) plane wave and local complex Taylor polynomial expansions for much larger
leaf boxes in the numerical discretization. However there are several numerical challenges: 
It becomes possible to have new poles (roots of $z+\mathrm{i} s(z)-w =0$) moving inside the 
contour for $I_1$ or $I_2$; more polynomial expansion terms become necessary for a prescribed
accuracy requirement; computing the local expansion coefficients may require new contours 
to avoid any numerical difficulties; and the discretization scheme to generate the FMM hierarchical 
tree structure also requires further study for optimal performance of the algorithm. 
These new challenges are being studied.

\subsection{Laplace Single Layer Potential: Integration by Parts}
\label{sec:LaplaceSingleLayer}
To derive the representation for the single layer potential
\begin{equation}
\begin{array}{rl}
SLP(w) = & \Re (\frac{1}{2\pi} \int\limits_{-1}^1 \log |w-z| \rho (z) \mathrm{d} z) \\
  = &  \Re (\frac{1}{2\pi} \int\limits_{-1}^1 \frac{1}{2} \log (
        (x-\tilde{x})^2+ (y-s (\tilde{x}))^2) \rho (\tilde{x})\left|
        1+\mathrm{i} s' (\tilde{x}) \right| \mathrm{d} \tilde{x})\\
  =&  \frac{1}{2 \pi} \int\limits_{-1}^1 \frac{1}{2}\log (
     (x-\tilde{x})^2+ (y-s (\tilde{x}))^2 ) \tilde{\rho}
     (\tilde{x})\mathrm{d} \tilde{x}\\
\end{array} 
\end{equation}
where $w=x+\mathrm{i} y$, $z=\tilde{x}+\mathrm{i} s(\tx)$, and 
$\tilde{\rho} (\tilde{x})=  \rho (\tilde{x})\left| 1+\mathrm{i} s' (\tilde{x}) \right|.$
We first apply the Fourier extension technique to represent the real function
$\tilde{\rho} (\tilde{x})$ as
$\tilde{\rho} (\tilde{x})= \sum\limits_{p=-P}^P  \omega_p e^{\mathrm{i} p \tilde{x}},$
and define $$f (\tilde{x})= \sum\limits_{p=-P, p \neq 0}^P \frac{ \omega_p}{\mathrm{i} p}
e^{\mathrm{i} p \tilde{x}} + \omega_0 \tilde{x},$$
which is a particular anti-derivative of $\tilde{\rho} (\tilde{x})$ as 
$f' (\tilde{x})=\tilde{\rho}(\tilde{x})$. Using $f$ and integration by part, we have
\begin{equation}
\begin{array}{rl}
SLP (w)= &  \frac{1}{2\pi} \int\limits_{-1}^1 \frac{1}{2} \log
           ((x-\tilde{x})^2+ (y-s (\tilde{x}))^2)\mathrm{d} f (\tilde{x}) \\
=  &\frac{1}{2\pi} \frac{1}{2}  \log
           ((x-\tilde{x})^2+ (y-s (\tilde{x}))^2) f
     (\tilde{x})|_{-1}^{1} \\
& +  \frac{1}{2 \pi} \int\limits_{-1}^1 f (\tilde{x})
  \frac{(x-\tilde{x})+ (y-s (\tilde{x}))s'
  (\tilde{x})}{(x-\tilde{x})^2+ (y-s (\tilde{x}))^2} \mathrm{d}
  \tilde{x}\\
=& I_1+I_2.\\ 
\end{array}
\end{equation} 
For the $I_1$ term, as both end points are well-separated from the leaf box 
containing $w=x+\mathrm{i} y$, a local expansion can be derived.  
For $I_2$, simple algebra shows that 
\begin{equation}
\begin{array}{rl}
I_2= &  \frac{1}{2 \pi} \int\limits_{-1}^1 f (\tilde{x})
  \frac{(x-\tilde{x})}{(x-\tilde{x})^2+ (y-s (\tilde{x}))^2} \mathrm{d}
  \tilde{x} 
   +\frac{1}{2 \pi} \int\limits_{-1}^1 f (\tilde{x})
  \frac{ (y-s (\tilde{x}))s'
  (\tilde{x})}{(x-\tilde{x})^2+ (y-s (\tilde{x}))^2} \mathrm{d}
  \tilde{x} \\
=& \frac{1}{2\pi}\Re \int\limits_{-1}^1  (\frac{1}{w-
   (\tilde{x}+\mathrm{i} s (\tilde{x}))}) f (\tilde{x}) \mathrm{d}
   \tilde{x}
- \frac{1}{2 \pi}\Im \int\limits_{-1}^1  (\frac{1}{w- (\tilde{x}-s
  (\tilde{x}))})s' (\tilde{x})f (\tilde{x})\mathrm{d} \tilde{x}. \\
\end{array}
\end{equation}
Both terms are in the form of Eq.~(\ref{eq:dlp}) with different density functions, therefore
results from previous section for double layer potentials can be applied directly. We 
skip the details.

\subsection{Comparing Plane Wave Expansion with Complex Polynomial Expansion}
The new QB2X technique uses both the local complex Taylor polynomial
expansion and (the slightly modified) plane wave expansion. Different error sources
in the QB2X can be easily analyzed: when the density function for $-1 \leq \tx \leq 1$ is 
approximated by an orthogonal polynomial, the truncation error is a well-studied
topic; the error from the Fourier approximation is controlled by the Fourier extension
precomputation process which generates the linear mapping from the polynomial basis 
to Fourier basis; the error from the truncated local complex polynomial approximation 
of the far-field density contribution on the well-separated boundary segments or on 
the contour follows the standard FMM analysis when a proper contour is chosen; 
the (slightly modified) plane wave in the representation is derived analytically 
using the Residue Theorem; finally finding the roots of a polynomial is a well 
studied topic.

\begin{figure}
  \centering
    \includegraphics[width=0.46\linewidth]{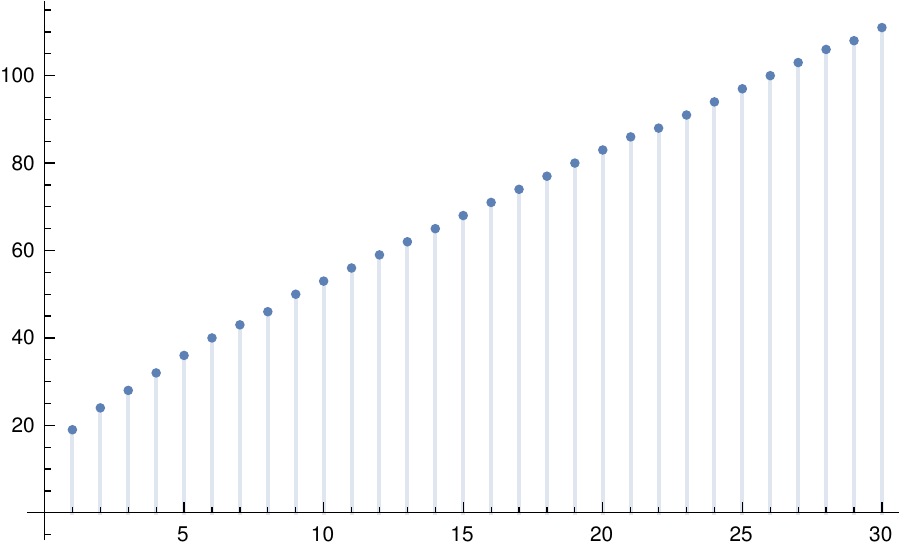}
  \caption[]
  {Estimated $N$ (y-axis) for different wave number $p$ (x-axis).  }
  \label{fig:kN}
\end{figure}

The QB2X also provides an alternative approach to analyze the error in the classical 
QBX method (see \cite{epstein2013convergence} for existing results), by studying 
the truncation error when a plane wave term $e^{ipw}$ is re-expanded as 
a local Taylor polynomial expansion 
\begin{equation}
  e^{\mathrm{i} p (w-w_0)}=\sum\limits_{k=0}^{\infty} \frac{(\mathrm{i} p)^k}{k!} (w-w_0)^{k}
\end{equation}
where $ w_0 $ is the expansion center. 
We assume $|w- w_0| \leq 1$, and study how many terms are required in the polynomial
expansion in order to achieve machine precision accuracy, i.e., we need to find
$N$ such that 
$|\sum\limits_{k=N}^{\infty} \frac{(\mathrm{i} p)^k}{k!} (w-w_{0})^{k}|\leq  10^{-16}.$
Using the incomplete gamma function $\Gamma
(N,p)=\int\limits_p^{\infty} e^{-x} x^{N-1}\mathrm{d} x$ and gamma
function $\Gamma(N)=\Gamma(N,0)$, we have
$$|\sum\limits_{k=N}^{\infty} \frac{(\mathrm{i}
    p)^k}{k!} (w-w_{0})^{k}|\leq \sum\limits_{k=N}^{\infty} \frac{p^k}{k!}= 
  \frac{e^p (\Gamma (N)-\Gamma (N,p))}{\Gamma (N)}.$$ In Fig.~\ref{fig:kN}, we numerically
solve the inequality $\frac{e^p (\Gamma (N)-\Gamma (N,p))}{\Gamma (N)} \leq 10^{-16}$ to
get an estimate of $N$. When $p=30$ is used in the plane wave expansion, to
achieve machine precision the estimated $N$ is about $111$, and this number becomes
larger if the nonlinearity due to the boundary geometry $s(\tilde{x})$ is included. 
Therefore, introducing two different bases in the representation improves the 
efficiency when approximating the layer potentials as we may gain accuracy from the plane wave approximation that would take many terms from a complex Taylor polynomial expansion alone.  

\section{Numerical Experiments}
\label{sec:num}
We present preliminary numerical results to validate the analytical formulas and demonstrate
the achieved accuracy for different $K$ (for the polynomial expansion) and $P$ 
(for the plane wave expansion) values.

\subsection{Double Layer Potential}
\begin{figure}[htbp]
  \centering
  \subfigure[$\rho (x)=\cos (x)$, $ P = 1$]{
    \includegraphics[width=0.46\linewidth]{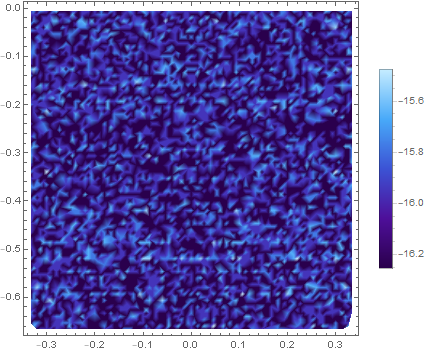}
  }
  \subfigure[ $\rho (x)=e^{\cos (x)}$, $ P = 20$]{
    \includegraphics[width=0.46\linewidth]{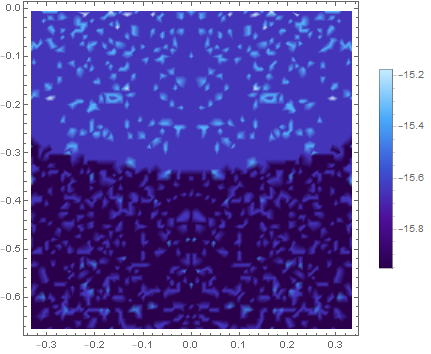}
  }
  \subfigure[ $\rho (x)=\frac{1}{4} \left(2 x^2+2 x+3\right)$, $P = 30$]{
    \includegraphics[width=0.46\linewidth]{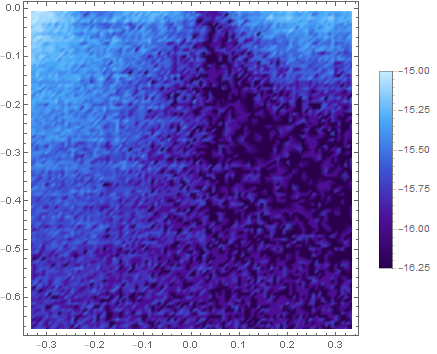}
  }
  \subfigure[ $\rho (x)=\frac{1}{8} \left(4 x^3+4 x^{2}+x+6\right)$, $P = 30$]{
    \includegraphics[width=0.46\linewidth]{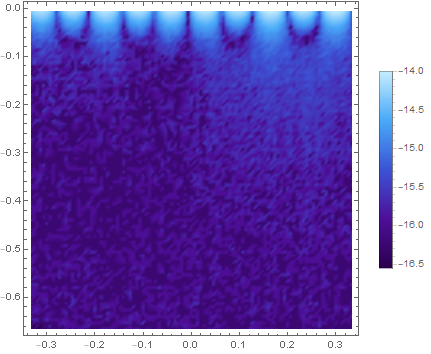}
  }
  \caption[]
  {Approximation errors of QB2X representations for double layer potentials with $ K = 40$ for 
    different density functions.  The plot legends are $\log_{10} (Error)$.}
  \label{fig:DLPerrStraight}
\end{figure}

We first consider the straight line segment connecting $(-1,0)$ and $(1,0)$ when $s(\tilde{x})=0$.   
The double layer potential of interest is then given by
\begin{equation}
DLP (w)= \frac{1}{2\pi} \Im \left( \int\limits_{-1}^1  \frac{ -\rho (\tilde{x})}{\tilde{x}-w} \mathrm{d} \tilde{x} \right).
\end{equation}
In existing implementations which combine QBX with FMM, the density function $\rho(\tilde{x})$ 
is often approximated by an orthogonal polynomial expansion, e.g., the Chebyshev polynomial expansion 
$\rho (\tilde{x})=\sum\limits_{n=0}^{N} c_n T_n (\tilde{x})$ where $T_n (\tilde{x})$ is the $n_{th}$ 
Chebyshev basis polynomial given by $T_n (\cos \theta) =\cos(n \theta).$ The first three basis polynomials 
in terms of $\tx=\cos(\theta)$ are explicitly given by $T_0(\tx)=1$, $T_1(\tx)=\tx$, and $T_2(\tx)=2\tx^2-1$. 
In the first numerical test, we choose $4$ different density $\rho$ functions: (a) $\rho (\tx)=\cos (\tx) $; 
(b) $\rho (\tx) = e^{\cos (\tx)}$; (c) $\rho (\tx)=T_0 (\tx)+\frac{1}{2}T_1(\tx)+ \frac{1}{4}T_2 (\tx)
=\frac{1}{4} \left(2 \tx^2+2 \tx+3\right)$; and (d) $\rho (\tx)=\frac{1}{8} \left(4 \tx^3+4 \tx^{2}+ \tx+6\right)$.
For (a), as it is already in the form of an exponential expansion, so $P=1$. We use $P=20$ in the Fourier expansion
for (b) to guarantee machine precision accuracy. For (c) and (d), we apply the precomputed mapping from the 
polynomial basis to the Fourier basis to compute the Fourier extensions with $P=30$. The approximation errors
are also around machine precision. We assume the target point 
$w \in [-\frac{1}{3}, \frac{1}{3}]\times [-\frac{2}{3},0]$ and the center of the leaf box is 
given by $w_0=(0, \frac{-1}{3})$. For the $I_1$ and $I_2$ terms given explicitly in Eq.~(\ref{eq:I1I2}),
we choose the contours in Fig.~\ref{fig:contours}, where $L \to \infty$ is used.
The $r_{max}$ values are $r_{max}=0.354$ for the upper contour and $r_{max}=0.471$ for the lower contour,
respectively. We choose $K=40$ which guarantees at least $13$-digits accuracy in the complex local 
Taylor polynomial expansion using standard FMM error analysis. The approximation errors are shown in the 
Fig.~\ref{fig:DLPerrStraight}. For all cases, the QB2X representations achieve $14$-digits accuracy.

\begin{figure}
  \centering
  \subfigure[$K=9$]{
    \includegraphics[width=0.46\linewidth]{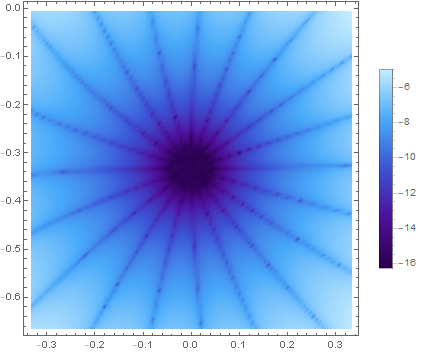}
  }
  \subfigure[$K=18$]{
    \includegraphics[width=0.46\linewidth]{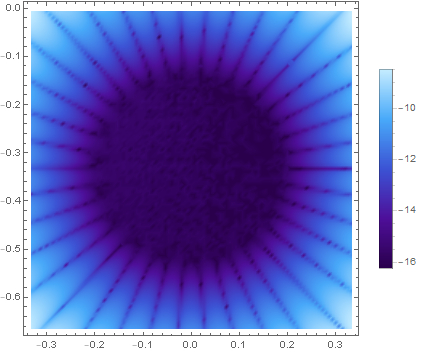}
  }
  \subfigure[$K=27$]{
    \includegraphics[width=0.46\linewidth]{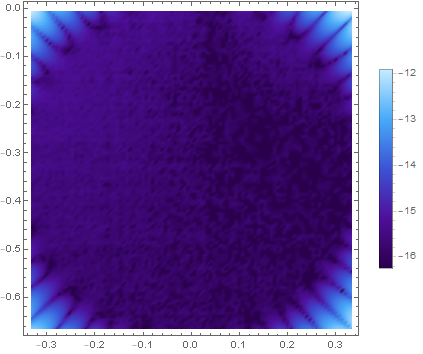}
  }
  \subfigure[ $K=36$]{
    \includegraphics[width=0.46\linewidth]{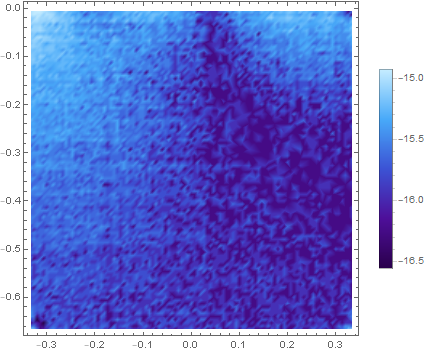}
  }
  \caption[]
  {Approximation errors of the QB2X double layer potential representations with 
    $\rho (\tx)=\frac{1}{4} \left(2 \tx^2+2 \tx+3\right)$ for different $K$ values. 
	$P= 30$. The plot legends are $\log_{10} (Error)$.}
  \label{fig:DLPerrStraightK}
\end{figure}

In classical FMM error analysis, when $K=9$, the complex local Taylor polynomial expansion is guaranteed to
achieve $3$-digits accuracy, and the accuracies increase to $6$, $9$, and $12$ digits when
$K=18$, $27$, and $36$, respectively. The same error estimates can be derived using the $r_{max}$ values in this
example. In Fig.~\ref{fig:DLPerrStraightK}, we show the approximation error for
different number of expansion terms $K$ for case (c) when $\rho (\tx) =\frac{1}{4} \left(2 \tx^2+2 \tx+3\right)$. 
In the experiment, we fix $P=30$ so the error from the Fourier extension is within machine precision. 
For all tested $K$ values, the errors are less than the error bound estimates derived using $r_{max}$.

\begin{figure}
  \centering
  \subfigure[$s(\tx)=-\frac{1}{10} \tx^{2}, K=40$]{
   \includegraphics[width=0.46\linewidth]{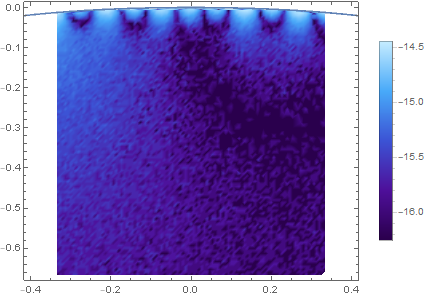}
  }
  \subfigure[$s(\tx)=-\frac{1}{10} \tx^{2}-\frac{1}{10} \tx^{4}, K=50$]{
    \includegraphics[width=0.46\linewidth]{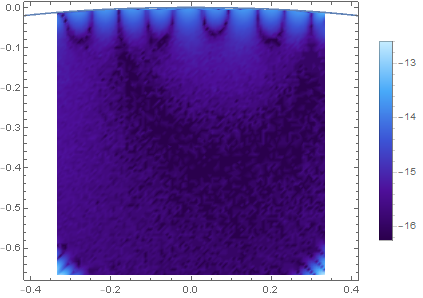}
  }
  \caption[]
  {Approximation errors of double layer potential for different boundary curves. 
  $P = 30$, the density function is 
  $\rho (\tx)=\frac{1}{4} \left(2 \tx^2+2 \tx+3\right)$, and the boundary curves are 
	$(\tx,-\frac{1}{10} \tx^{2})$ (left) and $( \tx,-\frac{1}{10} \tx^{2}-\frac{1}{10} \tx^{4})$ (right), 
	respectively, where $ \tx \in [-\frac{1}{3},\frac{1}{3}]$. The plot legend is $\log_{10} (Error)$.}
  \label{fig:DLPerrCurve}
\end{figure}

Next we study two curved boundaries defined by $s(\tx)=-\frac{\tx^2}{10}$ and 
$s(\tx)= - \frac{\tx^2}{10} -\frac{\tx^4}{10} $, respectively. We choose 
$\rho (\tx)=\frac{1}{4} \left(2 \tx^2+2 \tx+3\right)$ for both curves and consider
the leaf target box $w \in \left\{ (x,y)| \frac{-1}{3}<x<\frac{1}{3}, -\frac{2}{3}<y< s(x) \right\}$ 
with center $w_0= (0,-\frac{1}{3})$. The double layer potential is 
\begin{equation}
	DLP (w)=  \frac{1}{2\pi} \Re\int\limits_{-1}^1  \frac{ s'(\tilde{x}) \rho (\tilde{x})}{w- (\tilde{x}+
	\mathrm{i} s(\tilde{x}))} \mathrm{d} \tilde{x}+
           \frac{1}{2\pi} \Im\int\limits_{-1}^1  \frac{ \rho (\tilde{x})}{w-
(\tilde{x}+\mathrm{i} s(\tilde{x}))} \mathrm{d} \tilde{x}.
\end{equation}
Let $f_1 (\tx)=s'(\tx) \rho (\tx)$ and $f_2 (\tx)=\rho(\tx)$, we computed the QB2X representations of 
$DLP (w)$ for both boundary curves using the method in Sec.~\ref{sec:doubleCurve}. 
We present the errors in Fig.~\ref{fig:DLPerrCurve}. On the left plot, 
the $r_{max}$ values for the boundary curve $s(\tx)= -\frac{\tx^2}{10}$ are 
$r_{max}=0.329$ for the upper contour and $r_{max}=0.506$ for the lower contour, 
respectively, therefore $K=40$ will guarantee $12$ digits accuracy. On the right, 
the $r_{max}$ values for $s(\tx)= - \frac{\tx^2}{10} -\frac{\tx^4}{10} $ are 
$r_{max}=0.354$ for the upper contour and $r_{max}=0.622$ for the lower contour, 
respectively, therefore $K=40$ only givens numerical result with about $10$-digits 
accuracy, and a larger number $K=50$ is required for $12$ digits accuracy.

\subsection{Single Layer Potential}
\begin{figure}
  \centering
  \subfigure[$\rho(\tx)=\cos (\tx)$, $ P = 1$]{
    \includegraphics[width=0.46\linewidth]{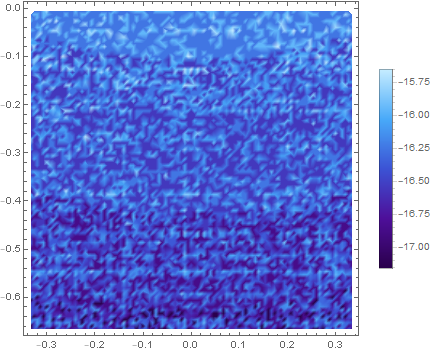}
  }
  \subfigure[ $\rho (\tx)=e^{\cos (\tx)}$, $ P = 20$]{
    \includegraphics[width=0.46\linewidth]{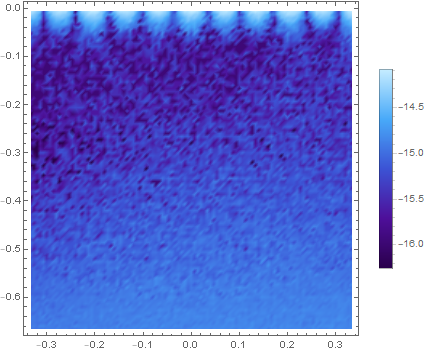}
  }
  \subfigure[ $\rho (\tx)=\frac{1}{4} \left(2 \tx^2+2 \tx+3\right)$, $P = 30$]{
    \includegraphics[width=0.46\linewidth]{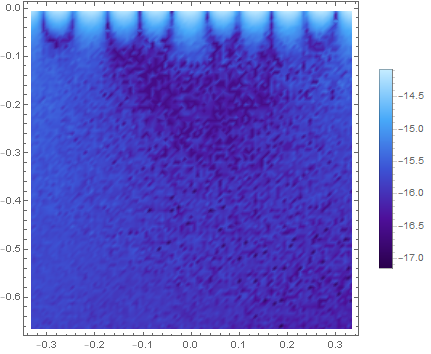}
  }
  \subfigure[ $\rho (\tx)=\frac{1}{8} \left(4 \tx^3+4 \tx^{2}+\tx+6\right)$, $P = 30$]{
    \includegraphics[width=0.46\linewidth]{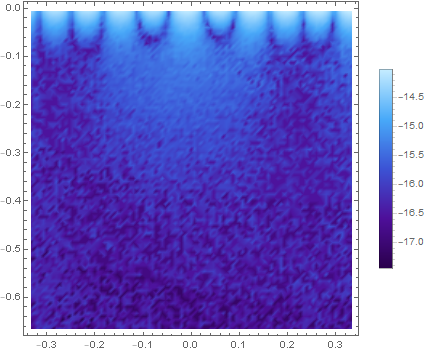}
  }
  \caption[]
  {Approximation errors of single layer potential for $ K = 40$ and
    different density functions. The boundary is a line segment
    connecting $(-1,0)$ and $(1,0)$. The plot legends are $\log_{10} (Error)$.}
  \label{fig:SLPerrStraight}
\end{figure}

Next we consider the QB2X representations for the single layer potentials.
We firstly consider the straight line case for different density functions defined 
on the line segment connecting $(-1,0)$ and $(1,0)$. We use the same density
functions as those in the double layer case, and the approximation errors are 
shown in Fig.~\ref{fig:SLPerrStraight} when $K=40$. For all cases, we
achieve $13$-digits accuracy.

\begin{figure}
  \centering
  \subfigure[$K=9$]{
    \includegraphics[width=0.46\linewidth]{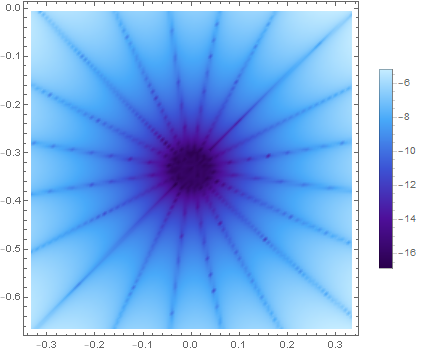}
  }
  \subfigure[$K=18$]{
    \includegraphics[width=0.46\linewidth]{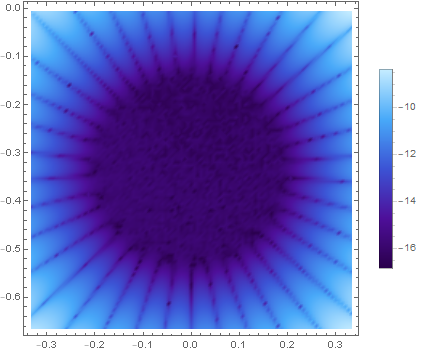}
  }
  \subfigure[$K=27$]{
    \includegraphics[width=0.46\linewidth]{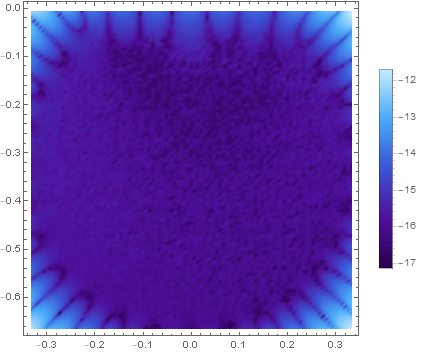}
  }
  \subfigure[ $K=36$]{
    \includegraphics[width=0.46\linewidth]{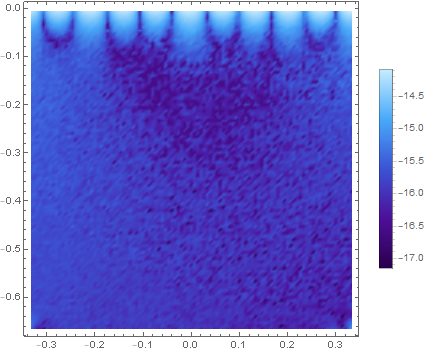}
  }
  \caption[]
  {Approximation errors of single layer potential for $ P= 30$ and
     density function $\rho (\tx)=\frac{1}{4} \left(2 \tx^2+2 \tx+3\right)$. The boundary is a line segment
    connecting $(-1,0)$ and $(1,0)$. The plot legends are $\log_{10} (Error)$.}
  \label{fig:SLPerrStraightK}
\end{figure}

To demonstrate the error dependency on the number of local Taylor polynomial 
expansion terms $K$, in Fig.~\ref{fig:SLPerrStraightK}, we plot the errors for different 
$K$ values when  $\rho (\tx)=\frac{1}{4} \left(2 \tx^2+2 \tx+3\right)$. 
Similar to the double layer case, for all tested $K$ values, the errors are smaller 
than the estimated bounds.

\begin{figure}
  \centering
  \subfigure[$K=18$]{
    \includegraphics[width=0.46\linewidth]{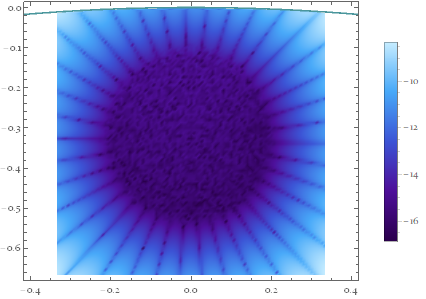}
  }
  \subfigure[$K=36$]{
    \includegraphics[width=0.46\linewidth]{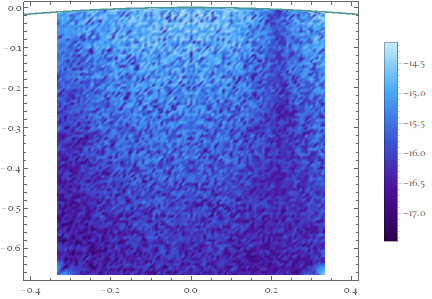}
  }
  \caption[]
  {Approximation errors of single layer potential for different $K$, $P =
    30$ and density functions  $\rho (\tx)=\frac{1}{4} \left(2 \tx^2+2
  \tx+3\right)$. The boundary is $(\tx,-\frac{1}{10} \tx^{2})$. 
	The plot legend is $\log_{10} (Error)$.}
  \label{fig:SLPerrCurve}
\end{figure}

Finally we consider the single layer potential on a curved boundary $(\tx,-\frac{1}{10} \tx^{2})$ when  
$\rho (\tx)=\frac{1}{4} \left(2 \tx^2+2 \tx+3\right)$ for different $K$ values. The single layer potential 
is given by
\begin{equation}
SLP(w) = \frac{1}{4 \pi} \int\limits_{-1}^1 \log (
     (x-\tilde{x})^2+ (y+\frac{1}{10}\tilde{x}^{2})^2 ) \tilde{\rho} (\tilde{x})\mathrm{d} \tilde{x}
\end{equation}
where $\tilde{\rho} (\tilde{x})=\rho (\tx) \sqrt{1+\frac{\tilde{x}^{2}}{25}}$. In the left plot of 
Fig.~\ref{fig:SLPerrCurve}, we show the error distribution when $K=18$ terms are used in the
polynomial expansion. The computed numerical results achieve at least $6$-digits accuracy. 
In the right plot, we show the error when $K=36$, and the results have at least $12$-digits
accuracy. 

\section{Summary}
\label{sec:summary}
In this paper, we present a new quadrature by two expansions (QB2X) technique for the Laplace layer 
potentials in two dimensions. Both the local complex Taylor polynomial expansion and 
plane wave expansions are used in the new representation. Compared with the classical QBX,
the new QB2X representations allow easier error analysis. For a prescribed accuracy requirement,
the QB2X representations are valid in a much larger region when compared with classical QBX 
representations. The impacts of the boundary geometry also become explicit in the QB2X 
representations, providing a useful tool for PDE analysis.

The QB2X technique can be generalized to other types of equations (e.g., the Helmholtz
and Yukawa equations) in both two and three dimensions using the Green's Identities. 
In these cases, the local complex Taylor polynomial expansions become the well-know
partial wave expansions. The partial wave and plane wave basis functions form a frame,
and the combined QB2X representations should have improved accuracy, stability and 
efficiency properties and can be easily combined with existing fast multipole methods 
when solving boundary value elliptic PDE problems. Results along these directions will 
be presented in subsequent papers.

\section*{Acknowledgement} 
We thankfully acknowledge the generous support of the NSF grants DMS1821093 (J. Huang) and NSF CAREER Grant DMS-1352353 
 (J.L. Marzuola).

\bibliography{qbe}
\bibliographystyle{abbrv}
\end{document}